\title{Complex hyperbolic surfaces whose closed geodesics are all simple}

\documentclass{article}
\usepackage{graphicx,pstricks}
\usepackage[all]{xy}
\usepackage{subfigure}
\usepackage{rotating}
\usepackage{amsfonts}
\usepackage{amssymb}
\usepackage{amsmath}
\usepackage{hyperref}

\usepackage{epsfig}
\usepackage{psfrag}
\usepackage{enumerate}
\usepackage{array}
\usepackage{theorem}

\newcommand{\C}{\mathbb{C}}
\newcommand{\R}{\mathbb{R}}

\newcommand{\Z}{\mathbb{Z}}

\newcommand{\X}{{\bf X}}

\newcommand{\HCd}{{\rm H}_\C^2}
\newcommand{\HCn}{{\rm H}_\C^n}

\newtheorem{thm}{Theorem}[section]

\newtheorem{cor}[thm]{Corollary}

\newtheorem{prop}{Proposition}[section]

{\theorembodyfont{\rmfamily} \newtheorem{rmk}[thm]{Remark}}

\newcommand{\Pf}{{\em Proof}. }
\newcommand{\EPf}{\hfill$\Box$\vspace{.5cm}}
\numberwithin{equation}{section}

\author{Sami Douba, D.~B.~McReynolds, Julien Paupert, Jonathan Vittore}

\begin{document}
\maketitle
\begin{abstract} 
\noindent We show that all closed geodesics in any arithmetic complex hyperbolic surface  of second type are simple. These include the fake projective planes studied by Mumford, Klingler, Prasad--Yeung, and Cartwright--Steger.
\end{abstract}
\section{Introduction}

Closed geodesics have long played a central role in our understanding of locally symmetric spaces. In the lowest dimension, complete, finite-area hyperbolic surfaces (with the exception of the thrice-punctured sphere) are well known to have infinitely many simple and non-simple closed geodesics. Huber's  prime geodesic theorem~\cite{H} states that the number of closed geodesics of length at most $L$ is asymptotic to $e^L/L$ for any such surface (closed geodesics will be assumed primitive throughout). In contrast, Mirzakhani \cite{Mi} showed that for any complete, finite-area, oriented hyperbolic surface $X$ of genus $g$ with $n$ cusps such that $3g-3+n>0$, the number of simple closed geodesics of length at most $L$ is asymptotic to $c_XL^{6g-6+2n}$ where $c_X>0$ depends on $X$. Taken together, these two results demonstrate that, in dimension two, ``most" closed geodesics are not simple. 
 
In higher dimensions the situation is more complicated. On one hand, any finite-volume locally symmetric space containing a properly immersed, totally geodesic, negatively curved surface contains all the geodesics of that surface, and hence contains infinitely many non-simple closed geodesics. In particular, any noncompact arithmetic locally symmetric space of noncompact type possesses infinitely many non-simple closed geodesics since, by the rational version of the Jacobson--Morosov lemma, any such locally symmetric space contains a properly immersed, totally geodesic, negatively curved surface. On the other hand, one would expect that geodesics in a higher-dimensional manifold generically do not intersect (themselves or each other). In fact, using arithmetic techniques, Chinburg and Reid \cite{CR} produced  infinitely many (pairwise non-commensurable) closed  hyperbolic $3$--manifolds all of whose  closed geodesics are simple. Interestingly, it was an open problem until very recently  whether or not \emph{every} closed hyperbolic $3$--manifold contains infinitely many simple closed geodesics; see \cite{LM} and \cite{HY} and Problem~3.16 of Kirby's list. 
 
In this paper we prove a result in the spirit of the Chinburg--Reid theorem for complex hyperbolic surfaces. 

\begin{thm}\label{Main} 
All closed geodesics in an arithmetic complex hyperbolic surface of second type are simple.
\end{thm}

We refer the reader to Section~\ref{AlgPrelim} for the definition of an arithmetic surface of second type. We note that this class contains infinitely many pairwise non-commensurable manifolds, and that these are all compact. It follows from the classification of arithmetic lattices in $\mathrm{PU}(2,1)$ (see \cite{DWM}) that this class includes all compact arithmetic complex hyperbolic surfaces whose fundamental groups have adjoint trace field $\mathbb{Q}$. Moreover, it follows from work of Klingler~\cite{K}, Yeung~\cite{Y}, Prasad--Yeung~\cite{PY} and Cartwright--Steger~\cite{CS} that any fake projective plane is arithmetic of second type. Recall that a \emph{fake projective plane} is a smooth, compact complex surface that has the same Betti numbers as the complex projective plane (namely 1, 0, 1, 0, 1) but is not biholomorphic to the latter. Mumford~\cite{Mu} produced the first example of a fake projective plane and showed that only finitely many exist. Much later, building on the papers listed above, Cartwright--Steger \cite{CS} classified all fake projective planes, showing that there are precisely 100 up to biholomorphism (50 up to homeomorphism).

More broadly, among all arithmetic lattices in noncompact simple Lie groups, those in the outer forms of type $\mathrm{A}_n$ that arise from cyclic algebras with unitary involutions are perhaps the most mysterious. Many of the fundamental questions for lattices remain open in this setting. In the case of real rank one and in the lowest dimension, these correspond to the arithmetic lattices of second type in the Lie groups isogenous to ${\rm SU}(2,1)$. There are by now several instances in which these lattices behave dramatically differently from those of first type. By work of Shimura \cite{Shimura} (in the noncocompact case) and Kazhdan \cite{Kazhdan} (in the cocompact case), every arithmetic lattice of first type in ${\rm SU}(2,1)$ has a finite-index subgroup $\Delta$ with $H^1(\Delta,\mathbb{Q}) \ne 0$. In contrast, Rogawski \cite{Rogawski} proved that every congruence subgroup $\Delta$ of a principal arithmetic lattice of second type in $\mathrm{SU}(2,1)$ satisfies $H^1(\Delta,\mathbb{Q}) = 0$. Moreover, the congruence subgroup problem is open for second-type lattices and known to fail for first-type lattices. Klingler \cite{K} proved that the fundamental groups of the fake projective planes are $\mathrm{PGL}_3$ super-rigid, whereas this fails in general for first-type lattices. The second-named author proved in \cite{McReynoldsTorsion} that second-type lattices in $\mathrm{PU}(p-1,1)$ can only contain cyclic torsion subgroups when $p$ is prime, whereas every finite group arises as a torsion subgroup in a first-type lattice in $\mathrm{PU}(p-1,1)$ for some prime  $p$. The first-type lattices are also known to have a rich supply of both complex and totally real, totally geodesic surface subgroups~\cite{CS2, PP}. When $p$ is prime, the orbifolds associated to the second-type lattices in $\mathrm{SU}(p-1,1)$ are known to have no proper, properly immersed, totally geodesic subspaces apart from closed geodesics. The main result of this paper provides another dichotomy between first- and second-type lattices in ${\rm SU}(2,1)$ via simple geodesics.

As noted above, all first-type arithmetic complex hyperbolic surfaces contain many properly immersed, totally geodesic submanifolds of real dimension two. Combining this fact with our result yields the following dichotomy for arithmetic complex hyperbolic surfaces.

\begin{cor}
For every arithmetic complex hyperbolic surface $M$, precisely one of the following holds: 
\begin{itemize}
\item[(i)]
$M$ has infinitely many non-simple closed geodesics (in which case $M$ is of first type), or
\item[(ii)]
$M$ has only simple closed geodesics (in which case $M$ is of second type).
\end{itemize}     
\end{cor}

One can compare our setting with the only other negatively curved symmetric space of real dimension~$4$, namely real hyperbolic $4$--space ${\rm H}_\R^4$. As all  arithmetic lattices in $\mathrm{Isom}({\rm H}_\R^4)$ arise from bilinear forms, the associated quotients all have a rich supply of properly immersed, totally geodesic $k$--submanifolds for ${k=1,2,3}$. In particular, these quotients always contain infinitely many distinct commensurability classes of properly immersed, totally geodesic surfaces~\cite{Mac,Rei}. One then deduces from Theorem~\ref{Main} that the arithmetic complex hyperbolic surfaces of second type are precisely the arithmetic negatively curved locally symmetric manifolds of real dimension $4$ all of whose closed geodesics are simple.

The proof of Theorem~\ref{Main} is in two parts. First, we show that any loxodromic element in an arithmetic lattice of second type in ${\rm SU}(2,1)$ has infinite-order holonomy (Proposition~\ref{IrrationalAngle}). Second, we show, using explicit trace identities due to Will~\cite{W} relating the trace of a loxodromic element of the form $A^{-1}B$ with the eigenvalues and complex cross-ratios associated to the fixed points of loxodromics $A$ and $B$, that if a complex hyperbolic surface $M$ contains a non-simple closed geodesic, then $\pi_1(M)$ must contain a loxodromic element with trivial holonomy (Corollary~\ref{NonSimpleHolonomy}). The result follows.

\subsection*{Acknowledgments} We thank Alan Reid and Matthew Stover for helpful discussions and comments on a previous draft of this note.

The second-named author used ChatGPT 5.6 Pro during the post-production stage to fix grammar, spelling, and typesetting issues. The authors take full intellectual responsibility for the content of this paper, which, mathematical or otherwise, was produced by the authors. 

\section{Geometric preliminaries}

\subsection{Complex hyperbolic space and isometries}

The standard reference for complex hyperbolic geometry is \cite{G}. For the reader's convenience, we include a brief summary of key definitions and facts. We will only consider the case of dimension $n=2$ in this paper, but the general setup is identical for higher dimensions, so we state it for all $n \geqslant 1$. 
Consider $\C^{n,1}$, the vector space $\C^{n+1}$ endowed with a Hermitian form $\langle \cdot \, , \cdot \rangle$ of signature $(n,1)$. Let $V^-=\left\lbrace Z \in \C^{n,1} | \langle Z , Z \rangle <0 \right\rbrace$.
Let $\pi\colon \C^{n+1}-\{0\} \longrightarrow \C{\rm P}^n$ denote the projectivization map. Define ${\rm H}_\C^n$ to be $\pi(V^-) \subset \C{\rm P}^n$, endowed with the Bergman metric $d$ given by:

\begin{equation}\label{dist}
\cosh ^2 \frac{1}{2}d(\pi(X),\pi(Y)) = \frac{|\langle X, Y \rangle|^2}{\langle X, X \rangle  \langle Y, Y \rangle}
\end{equation}
From this formula, it is clear that ${\rm PU}(n,1)$ acts by isometries
on ${\rm H}_\C^n$ (where ${\rm U}(n,1)$ is the subgroup of ${\rm GL}(n+1,\C)$ preserving $\langle \cdot , \cdot \rangle$, and ${\rm PU}(n,1)$ is its image in ${\rm PGL}(n+1,\C)$).

\medskip

{\bf Fact:} ${\rm Isom}^0({\rm H}_\C^n)={\rm PU}(n,1)$, and ${\rm Isom}({\rm H}_\C^n)={\rm PU}(n,1) \ltimes \Z/2$ (adjoining a real reflection, see below).
\medskip

{\bf Classification:} Any $g \in {\rm PU}(n,1) \setminus\{ {\rm Id} \}$ is of one of the following
types:
\begin{itemize}
\item 
\emph{elliptic}: $g$ has a fixed point in ${\rm H}_\C^n$.
\item 
\emph{parabolic}: $g$ has (no fixed point in ${\rm H}_\C^n$ and) exactly one fixed point in $\partial_\infty{\rm H}_\C^n$.
\item 
\emph{loxodromic}: $g$ has (no fixed point in ${\rm H}_\C^n$ and) exactly two fixed points in $\partial_\infty{\rm H}_\C^n$.
\end{itemize}

Chen and Greenberg gave the following characterization of conjugacy classes in ${\rm U}(n,1)$ in terms of matrix representatives (Theorem 3.4.1 of \cite{CG}; see Section 2.4 of \cite{PW} for more details):

\begin{thm}[Chen--Greenberg]\label{conjclasses} 
\quad
\begin{itemize}
\item[(a)] 
Any elliptic element is semisimple, with eigenvalues of norm 1. Two elliptic elements are conjugate in ${\rm U}(n,1)$ if and only if they have the same eigenvalues and the same eigenvalue of negative type.
\item[(b)] 
Any loxodromic element is semisimple, with exactly $n-1$ eigenvalues of norm 1. Two loxodromic elements are conjugate in ${\rm U}(n,1)$ if and only if they have the same eigenvalues.
\item[(c)] 
Any parabolic element is not semisimple, and all its eigenvalues have norm 1. It has a unique decomposition $g=pe=ep$ with $p$ unipotent and $e$ elliptic. Two parabolic elements are conjugate in ${\rm U}(n,1)$ if and only if their unipotent and elliptic components are conjugate.
\item[(d)] 
There are two classes of unipotent elements for $n>1$, the vertical Heisenberg translations and the non-vertical Heisenberg translations.
\end{itemize}
\end{thm}

Loxodromic and parabolic isometries have nice matrix representatives in the \emph{Siegel model} of ${\rm H}^2_\C$, that is, the projective model described above using the Hermitian form given by the matrix:

\begin{equation}\label{SiegelForn}
H=
\begin{bmatrix}
0& 0 & 1 \\
0 & 1 & 0 \\
1 & 0 & 0
\end{bmatrix}
\end{equation}

For concreteness, we note that in this model any loxodromic isometry $g$ has a lift to ${\rm SU}(2,1)$ that is conjugate to a matrix of the form 

\begin{equation}\label{LoxoMatrix}
L_{(r,\theta)}=
\begin{bmatrix}
re^{-i\theta/3} & 0 & 0 \\
0 & e^{2i\theta/3} & 0 \\
0 & 0 & e^{-i\theta/3}/r
\end{bmatrix}
\end{equation}
for some $r >1$ and $\theta \in \R/2\pi \Z$. Note that while the lift of an isometry $g \in {\rm PU}(2,1)$ to ${\rm SU}(2,1)$ is only defined up to multiplication by a cube root of unity, the quantities $r $ and $\theta$ appearing above are independent of this choice and are in fact the geometric invariants of the loxodromic isometry $g$: the quantity $r$ is related to the translation length $l$ of $g$ (specifically, by $r+1/r=2\cosh l/2$), and $\theta$ is the \emph{holonomy angle} (or rotation angle) of $g$. We will also call $e^{i\theta} \in {\rm U}(1)$ the \emph{holonomy} of $g$.

\medskip

{\bf Definitions:} A {\em complex $k$--plane} is a projective $k$--dimensional subspace of $\C P^n$ intersecting $\pi(V^-)$ non-trivially (so, it is an isometrically embedded copy of ${\rm H}_\C^{k} \subset {\rm H}_\C^n$). Complex 1--planes are usually called {\em complex lines}. A {\em complex reflection} is an elliptic isometry $g\in {\rm PU}(n,1)$ whose fixed-point set is a complex $(n-1)$--plane.  
 
A real $k$--plane is the projective image of a totally real $(k+1)$--subspace $W$ of $\mathbb{C}^{n,1}$, i.e., a $(k+1)$--dimensional real vector subspace such that $\langle v,w\rangle\in\mathbb{R}$ for all $v,w\in W$ and such that $W$ intersects $\pi(V^-)$ nontrivially; hence its projective image is an isometrically embedded copy of $\mathbf{H}^k_{\mathbb{R}}\subset \mathbf{H}^n_{\mathbb{C}}$. We will usually call real 2--planes simply {\em real planes}, or {\em $\R$--planes}. Every real $n$--plane in ${\rm H}_\C^n$ is the fixed-point set of an isometry of order 2 called a {\em real reflection} or {\em $\R$--reflection}. The prototype of such an isometry is the map given in affine coordinates by $(z_1,...,z_n) \mapsto (\overline{z_1},...,\overline{z_n})$ (this is an isometry provided that the Hermitian form has real coefficients).

By Proposition~2.5.1 of \cite{CG}, the totally geodesic submanifolds of complex hyperbolic space ${\rm H}^n_\C$ are exactly the real and complex $k$--planes with $0\leqslant k \leqslant n$.

\subsection{Quadruples of boundary points and cross-ratios}

An interesting feature of complex hyperbolic space is that, unlike in real hyperbolic space, triples of boundary points generically do not belong to the boundary of a common totally geodesic plane. In fact, triples of boundary points in $\HCn$ are classified up to isometry by a real parameter, the \emph{Cartan angular invariant}, and are contained in the boundary of a copy of ${\rm H}^1_\C$ (resp. ${\rm H}^2_\R$) exactly when this invariant is maximal (resp. 0), see Section 7.1 of~\cite{G} for more details.  

The complex cross-ratio of an ordered quadruple of points of $\partial_\infty \HCn$ was introduced by Kor\'{a}nyi and Reimann \cite{KR} in a different form; the definition below is due to Goldman. Properties of the complex cross-ratio are studied in detail in Section 7.2 of \cite{G} and Section~4 of \cite{W} to which we refer the reader for more details.

Given an ordered quadruple of distinct points $(p_1,p_2,p_3,p_4)$ in $\partial_\infty {\rm H}_\C^n$, define:

\begin{equation}\label{KRdef} 
{\bf X}(p_1,p_2,p_3,p_4)=\frac{\langle P_3,P_1 \rangle \langle P_4,P_2 \rangle}{\langle P_4,P_1 \rangle \langle P_3,P_2 \rangle}
\end{equation}
where $P_i$ is any lift of $p_i$ in $\C^{n,1}$ (for $1\leqslant i \leqslant 4$). Then the nonzero complex number ${\bf X}(p_1,p_2,p_3,p_4)$ is independent of the choice of lifts and well-defined (no two linearly independent vectors of nonpositive norm in $\C^{n,1}$ are orthogonal).

Note that this cross-ratio is preserved by holomorphic isometries, and sent to its complex conjugate by antiholomorphic isometries:

\begin{prop}\label{invX} 
Let $p_1,p_2,p_3,p_4$ be 4 distinct points in $\partial_\infty {\rm H}_\C^n$.
\begin{itemize}
\item For $g \in {\rm PU}(n,1)$, we have $\X(g(p_1),g(p_2),g(p_3),g(p_4))=\X(p_1,p_2,p_3,p_4)$ . 
\item If $g$ is an antiholomorphic isometry, then $\X(g(p_1),g(p_2),g(p_3),g(p_4))=\overline{\X(p_1,p_2,p_3,p_4)}$. 
\end{itemize}
\end{prop}

Given a set of 4 distinct points $\{p_1,p_2,p_3,p_4\}$ in $\partial_\infty {\rm H}_\C^n$, there are a priori 24 different cross-ratios, according to different orderings of the points. However, one easily sees from the definition that most of these will be related by simple functions, and in fact all 24 cross-ratios are determined by the following 3: $\X(p_1,p_2,p_3,p_4)$, $\X(p_1,p_4,p_2,p_3)$ and $\X(p_1,p_3,p_4,p_2)$. The effect on the cross-ratio of some permutations of the 4 points is as follows:

\begin{prop} 
Let $p_1,p_2,p_3,p_4$ be 4 distinct points in $\partial_\infty {\rm H}_\C^n$. Then:
\begin{itemize}
\item $\X(p_2,p_1,p_4,p_3)=\X(p_1,p_2,p_3,p_4)$
\item $\X(p_3,p_4,p_1,p_2)=\overline{\X(p_1,p_2,p_3,p_4)}$
\item $\X(p_2,p_1,p_3,p_4)=\X(p_1,p_2,p_3,p_4)^{-1}$
\end{itemize}
\end{prop}

We will only use these particular identities, but note that among the 24 cross-ratios, there are only three after identifying cross-ratios related by equality, complex conjugation, and inversion. Falbel showed that these three cross-ratios (satisfying two real equations) determine the quadruple of distinct boundary points up to isometry; see Section~4.2.1 of \cite{W} for more details.  

\section{Arithmetic lattices of second type in $\mathrm{SU}(2,1)$}\label{AlgPrelim}

We refer the reader to \cite{Sch} and \cite{DWM} for background material on this section. The commensurability classes of arithmetic lattices of second type in $\mathrm{SU}(2,1)$ are in bijection with certain $F$--forms of $\mathrm{SL}_3$ where $F$ is a totally real number field. We fix distinct real places $\nu_1,\dots,\nu_d\colon F \to \mathbb{R}$, where $d = [F:\mathbb{Q}]$. Let $E/F$ be a totally imaginary quadratic extension with complex places $\omega_1,\dots,\omega_d\colon E \to \mathbb{C}$ such that $(\omega_j)_{| F} = \nu_j$. We denote the completion of $F$ (respectively,~$E$) with respect to a place $\nu$ (resp.,~$\omega$) by $F_\nu$ (resp.,~$E_\omega$).

Let $A/E$ be a central, simple $E$--division algebra of degree 3 that admits an involution $\tau$ of second kind. There are infinitely many nonisomorphic $E$--algebras that admit such an involution. Moreover, $\tau$ is uniquely determined by the signatures at all of the real places, up to automorphisms of $A$. Given $h \in A^\times = A - \{0\}$, we say that $h$ is \emph{$\tau$--Hermitian} if $h = \tau(h)$. We can endow $A$ with a new involution of second kind $\tau_h$ defined by $\tau_h(x) = h \tau(x) h^{-1}$ and define
\[ \mathrm{U}(A,h) = \{ x \in A~:~x\tau_h(x) = 1\}, \quad \mathrm{SU}(A,h) = \mathrm{U}(A,h) \cap A^1, \] 
where $A^1$ is the group of reduced norm-one elements. We will denote the associated $F$--algebraic group by $\mathbf{G} = \mathbf{G}_{A,h}$ and note the group of $F$--points is $\mathrm{SU}(A,h)$. For each real place $\nu_j$, we have the group ${\mathrm{SU}(A \otimes_E E_{\omega_j},h)} = \mathbf{G}(F_{\nu_j})$. It follows that $\mathrm{SU}(A \otimes_E E_{\omega_j},h)$ is isomorphic either to $\mathrm{SU}(2,1)$ or to $\mathrm{SU}(3)$. In order to produce arithmetic lattices in $\mathrm{SU}(2,1)$, we require that there exists a real place $\nu_{j_0}$ such that
\begin{equation}\label{E:LatticeCond}
\mathrm{SU}(A \otimes_E E_{\omega_{j_0}},h) \cong \mathrm{SU}(2,1), \quad \mathrm{SU}(A \otimes_E E_{\omega_j},h) \cong \mathrm{SU}(3) \text{ for all }\nu_j \ne \nu_{j_0}.
\end{equation}
We now fix an $F$--form $\mathbf{G}(A,h)$ which satisfies (\ref{E:LatticeCond}). To produce a lattice, we must make one final choice, namely, we pick a maximal $\mathcal{O}_E$--order $\mathcal{O}$ in $A$. The group $\mathbf{G}(\mathcal{O}_F) = \mathrm{SU}(\mathcal{O},h) < \mathrm{SU}(A \otimes_E E_{\omega_{j_0}},h) \cong \mathrm{SU}(2,1)$ is a cocompact arithmetic lattice by Borel--Harish-Chandra \cite{BHC}. Note that there are only finitely many maximal orders up to conjugation and any other choice of a maximal order will produce a commensurable lattice. We say that a lattice $\Gamma$ in $\mathrm{SU}(2,1)$ is \emph{arithmetic of second type} if $\Gamma$ is commensurable in the wide sense with a lattice given as above. Likewise, a complex hyperbolic surface $M$ is {\em arithmetic of second type} if $M=\HCd/\Gamma$ for some arithmetic lattice $\Gamma < {\rm SU}(2,1)$ of second type.

The following proposition is known to experts. It follows, for instance, from Corollary~4.2 of~\cite{GG}; see also the proof of Lemma~2.9 of~\cite{RS}.

\begin{prop}\label{IrrationalAngle} 
If $\Gamma$ is an arithmetic lattice of second type in ${\rm SU}(2,1)$ and $\gamma \in \Gamma$ is loxodromic, then the holonomy of $\gamma$ has infinite order.
\end{prop}

\Pf By assumption, we have that $\Gamma$ is commensurable with a group of the form ${\rm SU}(\mathcal{O},h)$ as in the previous discussion. Hence, there is some integer $N > 0$ such that $\gamma^N \in A$ for some division algebra $A \subset \mathrm{M}(3,\mathbb{C})$. If the holonomy of $\gamma$ were of finite order, then after replacing $N$ be a multiple, we would have that $\gamma^N$ has $1$ as a simple eigenvalue, so that $\gamma^N - \mathrm{Id}$ would be a noninvertible nonzero element of $A$. But the latter is impossible since~$A$ is a division algebra.
\EPf

\section{Complex hyperbolic surfaces with a non-simple closed geodesic}

We start with a manifold $M=\HCd/\Gamma$ with $\Gamma$ an arithmetic lattice of second type in ${\rm SU}(2,1)$. If $M$ has a non-simple closed geodesic, lifting that geodesic to the universal cover provides two intersecting geodesics $\gamma_A$ and~$\gamma_B$ in $\HCd$, where $\gamma_A, \gamma_B$ are the respective axes of loxodromic isometries $A,B \in \Gamma$ with $A$ and $B$ conjugate in $\Gamma$.

The main result of this section is the following.

\begin{thm}\label{RealTrace} 
Let $A,B,C \in {\rm SU}(2,1)$ with $B=C^{-1}AC$ such that $A$ and $B$ are loxodromic isometries with intersecting axes. Then ${\rm Tr} \, A^{-1}B$ is real; in particular, if $ A^{-1}B$ is loxodromic, then its holonomy is $\pm 1$.
\end{thm}

\Pf Let $p_A,q_A$ (respectively $p_B,q_B$) denote the attracting and repelling fixed points of $A$ (respectively $B$) on $\partial_\infty \HCd$. Denote as in \cite{W}:
$$
\omega_a=\X(p_A,q_A,p_B,q_B)  , \  \omega_b=\X(p_A,q_B,q_A,p_B)  , \  \omega_c=\X(p_A,p_B,q_B,q_A).
$$ 
Now, since the axes $(p_A,q_A)$ and $(p_B,q_B)$ intersect, there exists a holomorphic isometry $g \in {\rm PU}(2,1)$ exchanging $p_A$ and $q_A$ as well as $p_B$ and $q_B$ (namely, $g$ is the central symmetry centered at the intersection point). It follows that
$$
\omega_b=\X(p_A,q_B,q_A,p_B)=\X(g(p_A),g(q_B),g(q_A),g(p_B))=\X(q_A,p_B,p_A,q_B)=\overline{\X(p_A,q_B,q_A,p_B)}=\overline{\omega_b}.
$$
Therefore, $\omega_b \in \R$. Analogously, we have $\omega_c \in \R$ since
$$
\omega_c=\X(p_A,p_B,q_B,q_A)=\X(g(p_A),g(p_B),g(q_B),g(q_A))=\X(q_A,q_B,p_B,p_A)=\overline{\X(p_A,p_B,q_B,q_A)}=\overline{\omega_c}
$$
We now use Will's explicit trace formulas (Section 5.2 of \cite{W}) for which we recall some notation. The \emph{complex multiplier} of a loxodromic element of ${\rm SU}(2,1)$ is its unique eigenvalue of modulus less than 1. For any two nonzero complex numbers $x,y$, denote $h(x,y)=(x-\overline{x}/x)(y-\overline{y}/y)$. Will's trace formula (Equation~(5.8) of \cite{W}) says that if $A,B$ are  loxodromic isometries  with complex multipliers $\mu,\nu$, respectively, then:
$$
{\rm Tr} \, A^{-1}B=\omega_c h(\mu^{-1},\nu)+\overline{\omega_c}h(\overline{\mu},\overline{\nu}^{-1}) + \frac{1}{\omega_b}h(\mu^{-1},\overline{\nu}^{-1})+\frac{1}{\overline{\omega_b}}h(\overline{\mu},\nu) + \beta(\mu,\nu) 
$$
where
$$
\beta(\mu,\nu)=\frac{\mu}{\overline{\mu}}(\nu+\overline{\nu}^{-1}) +\frac{\overline{\nu}}{\nu}(\mu^{-1}+\overline{\mu})-\frac{\mu \overline{\nu}}{\overline{\mu}\nu}.
$$ 
In our case the situation is much simplified as $\mu=\nu$ (because $A$ and $B$ are conjugate). The result follows by noting that $h$ is symmetric and that, for any $y,z \in \C$, we have $h(\overline{y},\overline{z})=\overline{h(y,z)}$, $h(z,\overline{z}) \in \R$, and $\beta(z,z) \in \R$. Then the terms appearing in the formula for ${\rm Tr} \, A^{-1}B$ are real or pairwise conjugate, and therefore ${\rm Tr} \, A^{-1}B \in \R$.

The second part follows from Proposition~2.3 of \cite{PW} (see also Theorem 6.2.4 of \cite{G}), as $ A^{-1}B$ then has eigenvalues $\{ \lambda, 1, \lambda^{-1} \}$ for some $\lambda \in \R$ with $|\lambda|>1$. If $\lambda>1$, then the holonomy of $ A^{-1}B$ is 1, and if $\lambda<-1$, the holonomy is $-1$.  \EPf

\begin{rmk} In the case of real hyperbolic 3--space studied by Chinburg and Reid, the intersecting axes must lie in a totally geodesic plane, forcing the four endpoints to have a real cross-ratio. They then use this to show that ${\rm Tr} \, [A,B^{-1}] \in \R$. The situation here is different:
\begin{itemize}
\item[(a)] The intersecting axes generically do not lie in a totally geodesic subspace. Indeed, given any three distinct points $p_1,p_2,p_3 \in \partial_\infty {\rm H}^2_\C$, there exists $p_4 \in \partial_\infty {\rm H}^2_\C$ such that the geodesics $(p_1, p_2)$ and $(p_3, p_4)$ intersect: just take $p_4 = g(p_3)$, where $g$ is the central involution centered at any point of $(p_1, p_2)$.
\item[(b)] In fact, here the condition ${\rm Tr} \, [A,B^{-1}] \in \R$ will only be satisfied when the four points $p_A,q_A,p_B,q_B$ lie on the boundary of a totally geodesic real plane in $\HCd$. Indeed, by the main result of \cite{W}, if ${\rm Tr} [A,B^{-1}] \in \R$ then there exists a real reflection $\sigma$ exchanging $p_A$ and $q_A$ as well as $p_B$ and $q_B$. But then $\sigma g$ would be an antiholomorphic isometry fixing these four points, which can only happen if  the four points lie on the boundary of a totally geodesic real plane.
\end{itemize}
\end{rmk}

\begin{rmk}
The conclusion of Theorem~\ref{RealTrace} does not hold if ${\rm SU}(2,1)$ is replaced by ${\rm PSL}(2,\mathbb{C})$, even when $A$ and $B$ are chosen within a cocompact lattice in ${\rm PSL}(2,\mathbb{C}) \cong {\rm Isom}^+({\rm H}_\R^3)$. Indeed, let $K \subset \mathbb{C}$ be a complex cubic field, and $\mathcal{O}_K$ be the ring of integers of $K$. Then $\Gamma:= {\rm SO}(3,\mathcal{O}_K)$ is a cocompact arithmetic lattice in ${\rm SO}(3,\mathbb{C}) \cong {\rm PSL}(2,\mathbb{C})$ and every loxodromic element of $\Gamma$ has infinite-order holonomy (equivalently, there is no infinite-order element of $\Gamma$ all of whose eigenvalues are real); see Theorem~5.3.1 in~\cite{MR}. On the other hand, the two diagonal matrices ${\rm diag}(1, -1, -1), {\rm diag}(-1,-1,1) \in \Gamma$ are conjugate within $\Gamma$ (indeed, within ${\rm SO}(3,\mathbb{Z})$), and so their centralizers within~$\Gamma$, which are each virtually infinite cyclic, are also conjugate within~$\Gamma$. If $A$ is a sufficiently high power of a loxodromic element of $\Gamma$ centralizing ${\rm diag}(1, -1, -1)$, and $B$ is an element of~$\Gamma$ conjugate to $A$ in $\Gamma$ and centralizing ${\rm diag}(-1,-1,1)$, then $A^{-1}B$ is loxodromic and, like all other loxodromic elements of $\Gamma$, has infinite-order holonomy. However, the axis of $A$ in ${\rm H}_\R^3$ intersects that of $B$; indeed, the axis of $A$ (respectively, the axis of $B$) is the fixed-point set of ${\rm diag}(1, -1, -1)$ (resp., of ${\rm diag}(-1,-1,1)$), and these fixed-point sets intersect, as can be seen from the fact that $AB$ again has finite order.
\end{rmk}

\begin{cor}\label{NonSimpleHolonomy} 
Assume $M={\rm H}^2_\C/\Gamma$ is a compact manifold; that is, $\Gamma \simeq \pi_1(M)$ is a torsion-free cocompact lattice in ${\rm PU}(2,1)$. If $M$ contains a non-simple closed geodesic, then $\Gamma$ contains a loxodromic element with trivial holonomy. 
\end{cor}

\Pf As noted above, if $M$ has a non-simple closed geodesic, lifting that geodesic to the universal cover provides two intersecting geodesics $\gamma_A$ and $\gamma_B$ in $\HCd$, where $\gamma_A$ (resp. $\gamma_B$) is the axis of a loxodromic isometry $A$ (resp. $B$) in $\Gamma$,  with $A$ and $B$ conjugate in $\Gamma$ (but not equal). 
Since $M$ is a compact manifold, every non-trivial element of $\Gamma$ is loxodromic; hence $A^{-1}B$ is loxodromic. By Theorem~\ref{RealTrace} (strictly speaking, applied to lifts of $A,B$ to SU(2,1)), we have that $A^{-1}B$ has holonomy $\pm 1$, so that $(A^{-1}B)^2$ has trivial holonomy.
\EPf

\raggedright

\quad

\begin{flushleft}
\textsc{Sami Douba\\
   Mathematisches Institut der Universit\"at Bonn}\\
   
       \verb|douba@math.uni-bonn.de|

       \medskip

  \textsc{D. B. McReynolds\\
   Purdue University}\\
       \verb|dmcreyno@purdue.edu|

       \medskip
       
  \textsc{Julien Paupert, Jonathan Vittore\\
   Arizona State University}\\
       \verb|paupert@asu.edu, jvittore@asu.edu|
\end{flushleft}


\frenchspacing
\begin{thebibliography}{ZZ99}

\bibitem[B-HC]{BHC}
A.~Borel, Harish-Chandra, {\sl Arithmetic subgroups of algebraic groups}, Ann. of Math. \textbf{75} (1962), 485--535.

\bibitem[CG]{CG} 
S.~Chen, L.~Greenberg, {\sl Hyperbolic spaces}, in Contributions to Analysis. Academic Press, New York (1974), 49--87.

\bibitem[CR]{CR} 
T.~Chinburg, A.~W.~Reid, {\sl Closed hyperbolic 3--manifolds whose closed geodesics all are simple}, J. Diff. Geom. {\bf 38} (1993), 545--558.

\bibitem[CaS]{CS} 
D.~I.~Cartwright, T.~Steger, {\sl Enumeration of the 50 fake projective planes}, C. R. Math. Acad. Sci. Paris {\bf 348} (2010), no. 1--2, 11--13.

\bibitem[ChS]{CS2}
T.~Chinburg, M.~Stover, {\sl Geodesic curves on Shimura surfaces}, Topology Proc., {\bf 52} (2018), 113--121.

\bibitem[Gol]{G} 
W.~M.~Goldman, {\sl Complex Hyperbolic Geometry}, Oxford Mathematical Monographs. Oxford University Press (1999).

\bibitem[GG]{GG} 
S.~Garibaldi, P.~Gille, {\sl Algebraic groups with few subgroups}, J. Lond. Math. Soc. (2) {\bf 80} (2009), 405--430.

\bibitem[Hub]{H} 
H.~Huber, {\sl Zur analytischen Theorie hyperbolischer Raumformen und Bewegungsgruppen}, Math. Ann. {\bf 138} (1959), 1--26.

\bibitem[HY]{HY} 
X.~H.~Han, B.~Yang, {\sl Growth of Simple Closed Geodesics in closed Hyperbolic 3-Manifolds}, Preprint (2026). Available at: \href{https://arxiv.org/abs/2608.29855}{https://arxiv.org/abs/2608.29855}. 

\bibitem[Kaz]{Kazhdan}
D.~Kazhdan, {\sl Some applications of the Weil representation}, J. Analyse Mat. \textbf{32} (1977), 235–248.

\bibitem[Kli]{K} 
B.~Klingler, {\sl Sur la rigidit\'e de certains groupes fondamentaux, l’arithm\'eticit\'e des r\'eseaux hyperboliques complexes, et les ``faux plans projectifs"}, Invent. Math. {\bf 153} (2003), 105--143.

\bibitem[KR]{KR} 
A.~Kor\'{a}nyi, H.~M.~Reimann, {\sl The complex cross-ratio on the Heisenberg group}, Enseign. Math. {\bf 33} (1987), 291--300.

\bibitem[LM]{LM} 
Q.~Luo, V.~Markovic, {\sl Simple geodesics in closed hyperbolic manifolds}, Preprint (2026). Available at: \href{https://arxiv.org/abs/2608.29761}{https://arxiv.org/abs/2608.29761}.

\bibitem[Mac]{Mac}
C.~Maclachlan, {\sl Fuchsian subgroups of the groups $\mathrm{PSL}_2(O_d)$}, in {\sl Low-dimensional Topology and Kleinian groups}, ed. D. B. A. Epstein, L. M. S. Lecture Note Series {\bf 112} (1986), 305--311.

\bibitem[MR]{MR} 
C.~Maclachlan, A.~W.~Reid, {\sl The arithmetic of hyperbolic 3--manifolds}, vol. \textbf{219} GTM, Springer-Verlag, 2003.

\bibitem[McR]{McReynoldsTorsion}
D.~B.~McReynolds, {\sl Finite subgroups of arithmetic lattices in $\mathrm{U}(2,1)$}, Geom. Dedicata, \textbf{122} (2006), 135–144.

\bibitem[Mir]{Mi} 
M.~Mirzakhani, {\sl Growth of the number of simple closed geodesics on hyperbolic surfaces}, Ann. of Math. {\bf 168} (2008), 97--125.

\bibitem[Mum]{Mu} 
D.~Mumford, {\sl An algebraic surface with $K$ ample, $K^2 = 9,\  p_g = q = 0$}, Amer. J. Math. {\bf 101} (1979), 233--244.

\bibitem[PP]{PP}
J.~Parkkonen, F.~Paulin, {\sl A classification of $\mathbb{R}$-Fuchsian subgroups of Picard modular groups}, Conflu. Math. {\bf 10} (2018), no. 2, 75–92.

\bibitem[PW]{PW} 
J.~Paupert, P.~Will, {\sl Real reflections, commutators, and cross-ratios in complex hyperbolic space}, Groups Geom. Dyn. {\bf 11} (2017), 311--352.

\bibitem[PY]{PY} 
G.~Prasad, S.-K.~Yeung, {\sl Fake projective planes}, Invent. Math. {\bf 168} (2007), 321--370.

\bibitem[Rei]{Rei}
A.~W.~Reid, {\sl Totally geodesic surfaces in hyperbolic 3-manifolds}, Proc. Edinburgh Math. Soc. (2) 34, no. 1 (1991), 77--88.

\bibitem[RS]{RS}
A.~W.~Reid, M.~Stover, {\sl Complex hyperbolic 2-orbifolds with isolated singularities}, Glasg. Math. J. {\bf 68}(3) (2026), 460--469.


\bibitem[Rog]{Rogawski}
J.~Rogawski, {\sl Automorphic representations of unitary groups in three variables}, volume 123 of Annals of Mathematics Studies. Princeton University Press, 1990.

\bibitem[Sch]{Sch}
W.~Scharlau, {\sl Quadratic and Hermitian forms}, \textbf{270} Grundlehren der Mathematischen Wissenschaften, Springer-Verlag, 1985.

\bibitem[Shi]{Shimura}
G.~Shimura, {\sl Automorphic forms and the periods of abelian varieties},  J. Math.
Soc. Japan \textbf{31} (1979), 561--592.

\bibitem[Wil]{W} 
P.~Will, {\sl Traces, Cross-ratios and 2--generator Subgroups of ${\rm PU}(2,1)$}, Can. J. Math. {\bf 61} (2009), 1407--1436.

\bibitem[WM]{DWM}
D.~Witte Morris, {\sl Introduction to Arithmetic Groups}, Deductive Press, 2001.

\bibitem[Yeu]{Y}
S.-K.~Yeung, {\sl Integrality and arithmeticity of co-compact lattice corresponding to certain complex two-ball quotients of Picard number one}, Asian J. Math. {\bf 8} (2004), 107--129.

\end{thebibliography}
\end{document}